\documentclass[10pt,a4paper,leqno,twoside]{article}

\usepackage{amsmath,amsthm,amssymb}
\numberwithin{equation}{section}

\newtheorem{thm}{Theorem}[section]
\newtheorem{lem}[thm]{Lemma}

\theoremstyle{definition}
\newtheorem{example}[thm]{Example}

\theoremstyle{remark}
\newtheorem{remark}[thm]{Remark}

\newenvironment{ex}{\begin{example}}{\hfill{$\diamondsuit$}\end{example}}
\newenvironment{rem}{\begin{remark}}{\hfill{$\diamondsuit$}\end{remark}}
\newenvironment{ackn}{\subsection*{Acknowledgement}}{}
\newenvironment{keywords}{\smallskip\noindent{\bfseries Keywords:}}{}
\newenvironment{MSC}{\smallskip\noindent{\bfseries 2000 Mathematical Subject Classification:}}{}

\newcommand{\C}{\mathbb{C}}
\newcommand{\R}{\mathbb{R}}
\newcommand{\I}{\sqrt{-1}}

\newcommand{\kakko}[1]{\left(#1\right)}
\newcommand{\ckakko}[1]{\left\{#1\right\}}
\newcommand{\abs}[1]{\left|#1\right|}

\DeclareMathOperator{\tr}{tr}
\DeclareMathOperator{\Tr}{Tr}
\DeclareMathOperator{\diag}{diag}
\DeclareMathOperator{\sgn}{sgn}
\DeclareMathOperator{\Sym}{Sym}

\newcommand{\deq}{:=}
\newcommand{\eqsp}{\phantom{{}={}}}

\newcommand{\B}{\boldsymbol{B}}
\newcommand{\E}{\boldsymbol{E}}
\newcommand{\F}{\boldsymbol{F}}

\newcommand{\va}{\boldsymbol{a}}
\newcommand{\vb}{\boldsymbol{b}}
\newcommand{\vd}{\boldsymbol{d}}
\newcommand{\ve}{\boldsymbol{e}}
\newcommand{\vj}{\boldsymbol{j}}
\newcommand{\vo}{\boldsymbol{o}}
\newcommand{\vt}{\boldsymbol{t}}
\newcommand{\vu}{\boldsymbol{u}}
\newcommand{\vv}{\boldsymbol{v}}
\newcommand{\vx}{\boldsymbol{x}}

\newcommand{\symx}[1]{\operatorname{Sym}^{\times}_{#1}}

\newcommand{\esum}[1]{\left\|#1\right\|}%{\Sigma(#1)}
\newcommand{\card}[1]{\left|#1\right|}

\newcommand{\T}[1]{\mathcal{T}_{#1}}
\newcommand{\U}[1]{\boldsymbol{U}_{\!#1}}
\newcommand{\V}[1]{\boldsymbol{V}_{\!\!#1}}
\newcommand{\den}[1]{\mathcal{W}_{#1}}

\newcommand{\hA}[2]{\frac{1-u_{#1}^4u_{#2}^4}{(1-u_{#1}^4)(1-u_{#2}^4)}}
\newcommand{\hB}[1]{\frac{-u_{#1}^2}{1-u_{#1}^4}}

\newcommand{\e}{\varepsilon}
\newcommand{\conj}[1]{\overline{#1}}
\newcommand{\sym}[1]{\mathfrak{S}_{#1}}
\newcommand{\cyclic}[1]{\mathcal{C}_{#1}}

\title{\bfseries Special value formula\\ for the spectral zeta function\\ of the non-commutative harmonic oscillator}
\author{Kazufumi KIMOTO}
\date{March 30, 2009}

\begin{document}

\maketitle

\begin{abstract}
We calculate the special values of the spectral zeta function of the non-commutative harmonic oscillator,
and give a general formula for them as integrals of certain algebraic functions.
This is a generalization of the result by Ichinose-Wakayama (2005),
in which the first two special values are studied.

\begin{keywords}
Non-commutative harmonic oscillator, spectral zeta function, special values.
\end{keywords}

\begin{MSC}
%11M41, %Other Dirichlet series and zeta functions
14G10, %Zeta-functions and related questions
35P99. %None of the above, but in this section (Spectral theory and eigenvalue problems for partial)
\end{MSC}
\end{abstract}

%=========================================================================
\section{Introduction}
%=========================================================================

Let $Q$ be a differential operator defined by
\begin{align*}
Q=Q_{\alpha,\beta}
=\begin{pmatrix}\alpha & 0 \\ 0 & \beta\end{pmatrix}\kakko{-\frac12\frac{d^2}{dx^2}+\frac12x^2}
+\begin{pmatrix}0 & -1 \\ 1 & 0\end{pmatrix}\kakko{x\frac{d}{dx}+\frac12}.
\end{align*}
The system defined by $Q$ is called the \emph{non-commutative harmonic oscillator},
which was first introduced and investigated by Parmeggiani and Wakayama \cite{PW2001, PW2002, PW2002c} (see also \cite{P}).
We always assume that $\alpha,\beta>0$ and $\alpha\beta>1$.
Under this assumption, $Q$ becomes positive self-adjoint unbounded operator on $L^2(\R;\C^2)$
which has only a discrete spectrum.
Denote the eigenvalues of $Q$ by
\begin{align*}
0<\lambda_1\le\lambda_2\le\lambda_3\le\dots(\to\infty),
\end{align*}
and define the spectral zeta function of $Q$ by
\begin{align*}
\zeta_Q(s)=\sum_{n=1}^\infty\frac1{\lambda_n^s}.
\end{align*}
This series is absolutely convergent in the region $\Re s>1$, and defines a holomorphic function in $s$ there.
We call this function $\zeta_Q(s)$ the \emph{spectral zeta function} for the non-commutative harmonic oscillator $Q$,
which is introduced and studied by Ichinose and Wakayama \cite{IW2005a}.
The zeta function $\zeta_Q(s)$ is analytically continued to the whole complex plane
as a single-valued meromorphic function which is holomorphic except for the simple pole at $s=1$.
It is remarkable that $\zeta_Q(s)$ has `trivial zeros' at $s=0,-2,-4,\dots$.
When the two parameters $\alpha$ and $\beta$ are equal, then $\zeta_Q(s)$ essentially gives the Riemann zeta function $\zeta(s)$
(see Remark \ref{rem:degenerate_case}).

We are interested in the special values of $\zeta_Q(s)$,
that is, the values $\zeta_Q(s)$ at $s=2,3,4,\dots$.
The first two special values are calculated by Ichinose and Wakayama \cite{IW2005b} as
\begin{align*}
\zeta_Q(2)&=2\kakko{\frac{\alpha+\beta}{2\sqrt{\alpha\beta(\alpha\beta-1)}}}^{\!\!2}\\
&\quad\eqsp\times\Biggl(\zeta(2,1/2)+\kakko{\frac{\alpha-\beta}{\alpha+\beta}}^{\!\!2}
\int_{[0,1]^2}\frac{4du_1du_2}{\sqrt{(1-u_1^2u_2^2)^2+(1-u_1^4)(1-u_2^4)/(\alpha\beta-1)}}\Biggr),
\end{align*}
\begin{align*}
\zeta_Q(3)&=2\kakko{\frac{\alpha+\beta}{2\sqrt{\alpha\beta(\alpha\beta-1)}}}^{\!\!3}\\
&\quad\eqsp\times\Biggl(\zeta(3,1/2)+3\kakko{\frac{\alpha-\beta}{\alpha+\beta}}^{\!\!2}
\int_{[0,1]^3}\frac{8du_1du_2du_3}{\sqrt{(1-u_1^2u_2^2u_3^2)^2+(1-u_1^4)(1-u_2^4u_3^4)/(\alpha\beta-1)}}\Biggr),
\end{align*}
where $\zeta(s,x)=\sum_{n=0}^\infty(n+x)^{-s}$ is the Hurwitz zeta function.
They also give the contour integral expressions of these values using a solution of a certain singly confluent Heun differential equation.
Later, Ochiai \cite{O} gave an expression of $\zeta_Q(2)$ using the complete elliptic integral
or the hypergeometric function (see Remark \ref{rem:Apery-Heun}),
and the author and Wakayama \cite{KW2006} gave a formula for $\zeta_Q(3)$ similar to the Ochiai's one.

In this article, we calculate the special values of the spectral zeta function $\zeta_Q(n)$
of the non-commutative harmonic oscillator $Q$ for all positive integers $n>1$
and express them in terms of integrals of certain algebraic functions (Theorem \ref{thm:specialvalues}).
Roughly, they are of the form
\begin{align*}
\zeta_Q(n)=2\kakko{\frac{\alpha+\beta}{2\sqrt{\alpha\beta(\alpha\beta-1)}}}^{\!\!n}
\Biggl(\zeta(n,1/2)+R_n(\alpha,\beta)\Biggr),
\end{align*}
and $R_n(\alpha,\beta)$ vanishes when $\alpha=\beta$.
It is natural to expect that
the `remainder term' $R_n(\alpha,\beta)$ has interesting properties
since it reflects the `non-commutativity' of the operator $Q$ in some sense.
In fact, for instance,
an investigation of $R_2(\alpha,\beta)$ for $\zeta_Q(2)$ leads us to a connection with
various arithmetic objects such as Ap\'ery-number-like binomial sums,
elliptic curves and modular forms (see, e.g., \cite{KW2007}),
and we expect that we might have such a kind of connection in the case of $R_n(\alpha,\beta)$ in general.
See \cite{K2009} for the first attempt to study this generalized situation.

%=========================================================================
\section{Preliminaries}
%=========================================================================

Following to Ichinose-Wakayama \cite{IW2005b},
we explain how to calculate the special values of $\zeta_Q(s)$.
Put
\begin{align*}
B(x,y)&\deq A^{-1}\ckakko{\cos\frac{q(x^2-y^2)}2I+\sin\frac{q(x^2-y^2)}2J},\\
E(u,x,y)&\deq
\exp\kakko{-\frac{1+u^4}{2(1-u^4)}(x^2+y^2)+\frac{2u^2}{1-u^4}xy}\\
&=\exp\kakko{-\begin{pmatrix} x & y \end{pmatrix}
\begin{pmatrix}
\frac{1+u^4}{2(1-u^4)} & \frac{-u^2}{1-u^4} \\
\frac{-u^2}{1-u^4} & \frac{1+u^4}{2(1-u^4)}
\end{pmatrix}
\begin{pmatrix} x \\ y \end{pmatrix}},
\end{align*}
where
\begin{align*}
\varepsilon\deq\frac1{\sqrt{\alpha\beta}},\quad
q\deq\frac\varepsilon{\sqrt{1-\varepsilon^2}}=\frac1{\sqrt{\alpha\beta-1}}
\end{align*}
and
\begin{align*}
A\deq\begin{pmatrix}
\alpha & 0 \\ 0 & \beta
\end{pmatrix},\quad
I\deq\begin{pmatrix}
1 & 0 \\ 0 & 1
\end{pmatrix},\quad
J\deq\begin{pmatrix}
0 & -1 \\ 1 & 0
\end{pmatrix}.
\end{align*}
Notice that $0<\e<1$ and $q>0$.
Furthermore, put
\begin{align*}
\B(x_1,\dots,x_n)&\deq\tr\kakko{B(x_1,x_2)B(x_2,x_3)\dots B(x_n,x_1)},\\[0.5em]
\E(u_1,\dots,u_n;x_1,\dots,x_n)&\deq E(u_1,x_1,x_2)E(u_2,x_2,x_3)\dots E(u_n,x_n,x_1),\\
\F(u_1,\dots,u_n)&\deq\int_{\R^n}\E(u_1,\dots,u_n;x_1,\dots,x_n)\B(x_1,\dots,x_n)dx_1\dots dx_n,
\end{align*}
where $\tr$ represents the matrix trace.
By Ichinose and Wakayama \cite{IW2005b}, the integral kernel of the operator $Q^{-1}$ is given by
\begin{align*}
Q^{-1}(x,y)=\int_0^1 A^{-1/2}K(u,x,y)A^{-1/2}du,
\end{align*}
where
\begin{align*}
K(u,x,y)=\frac1{(1-\e^2)^{1/4}\sqrt\pi}\frac1{\sqrt{1-u^4}}E(u,(1-\e^2)^{1/4}x,(1-\e^2)^{1/4}y)\exp\kakko{\frac{\e(x^2-y^2)}2J}.
\end{align*}
Hence, for a positive integer $n$, we have
\begin{equation}
\begin{split}\label{eq:special_value_primitive}
\zeta_Q(n)&=\Tr Q^{-n}\\
&=\int_{[0,1]^n}\int_{\R^n}\tr\kakko{A^{-1/2}K(u_1,x_1,x_2)A^{-1/2}\dotsb A^{-1/2}K(u_n,x_n,x_1)A^{-1/2}}d\vx d\vu\\
&=\kakko{\frac2{\sqrt{\pi(1-\varepsilon^2)}}}^{\!\!n}
\int_{[0,1]^n}\F(u_1,\dots,u_n)\frac{d\vu}{\sqrt{\prod_{j=1}^n(1-u_j^4)}},
\end{split}
\end{equation}
where $d\vx=dx_1\dots dx_n$, $d\vu=du_1\dots du_n$ and $\Tr$ is the operator trace.
Thus, we have only to calculate $\F(u_1,\dots,u_n)$
to get the special values of the spectral zeta function $\zeta_Q(s)$.

%=========================================================================
\section{Lemmas}
%=========================================================================

%-------------------------------------------------------------------------
\subsection{Calculation of the integrand of $\F(u_1,\dots,u_n)$}
%-------------------------------------------------------------------------

The following lemma is crucial.
\begin{lem}\label{lem:expansion_of_B}
For any positive integer $n$, it holds that
\begin{align*}
&\eqsp\B(x_1,x_2,\dots,x_n)\\
&=2\kakko{\frac{\alpha+\beta}{2\alpha\beta}}^{\!\!n}
\ckakko{1+\sum_{0<2k\le n}\kakko{\frac{\alpha-\beta}{\alpha+\beta}}^{\!\!2k}%
\sum_{1\le j_1<j_2<\dots<j_{2k}\le n}
\cos\kakko{q\sum_{r=1}^{2k}(-1)^{r}x_{j_r}^2}}.
\end{align*}
\end{lem}

\begin{proof}
Let us put $i=\sqrt{-1}$, $a_{1}=\alpha^{-1}$, $a_{2}=\beta^{-1}$ and $t_j=qx_j^2/2$ $(j=1,2,\dots,n)$.
We see that
\begin{align*}
&\eqsp\B(x_1,x_2,\dots,x_n)\\
&=\sum_{k_1,k_2,\dots,k_n\in\{1,2\}}
a_{k_1}a_{k_2}\dots a_{k_n}
\prod_{m=1}^n \cos\!\kakko{t_m-t_{m+1}+\frac{k_{m+1}-k_m}2\pi}\\
&=\sum_{k_1,k_2,\dots,k_n\in\{1,2\}}
a_{k_1}a_{k_2}\dots a_{k_n}
\prod_{m=1}^n \frac{i^{k_{m+1}-k_m}e^{i(t_m-t_{m+1})}+i^{-(k_{m+1}-k_m)}e^{-i(t_m-t_{m+1})}}2\\
&=\frac1{2^n}\sum_{k_1,k_2,\dots,k_n\in\{1,2\}}
\sum_{l_1,l_2,\dots,l_n\in\{1,-1\}}
\prod_{m=1}^n a_{k_m}i^{l_m(k_{m+1}-k_m)}e^{il_m(t_m-t_{m+1})}\\
&=\frac1{2^n}\sum_{l_1,l_2,\dots,l_n\in\{1,-1\}}
\sum_{k_1,k_2,\dots,k_n\in\{1,2\}}
\prod_{m=1}^n a_{k_m}i^{k_m(l_{m-1}-l_m)}e^{it_m(l_m-l_{m-1})}\\
&=\frac1{2^n}\sum_{l_1,l_2,\dots,l_n\in\{1,-1\}}
\prod_{m=1}^n (a_1i^{l_{m-1}-l_m}+a_2i^{2(l_{m-1}-l_m)})
e^{i(l_m-l_{m-1})t_m},
\end{align*}
where we set $k_0=k_n$, $k_{n+1}=k_1$, $l_0=l_n$, $l_{n+1}=l_1$, $t_0=t_n$ and $t_{n+1}=t_1$.
Here we notice that
\begin{itemize}
\item[(i)] $i^{l_{m-1}-l_m}=-(-1)^{\delta_{l_m,l_{m-1}}}$,
\item[(ii)] $\card{\{m\in\{1,2,\dots,n\}\,;\,l_{m-1}\ne l_m\}}$ is even (remark that $l_0=l_n$),
\item[(iii)] if there exist $j_1,\dots,j_{2k}\in\{1,2,\dots,n\}$ such that $j_1<\dots<j_{2k}$,
$l_{j_r-1}\ne l_{j_r}$ for $r=1,2,\dots,2k$ and $l_{m-1}=l_m$ for $m\in\{1,2,\dots,n\}\setminus\{j_1,\dots,j_{2k}\}$,
then $\sum_{m=1}^n (l_m-l_{m-1})t_m=2l_{j_1}\sum_{r=1}^{2k}(-1)^rt_{j_r}$.
\end{itemize}
Thus it follows that
\begin{align*}
&\eqsp\B(x_1,x_2,\dots,x_n)\\
&=\frac1{2^n}\sum_{l\in\{1,-1\}}\kakko{1+%
\sum_{0<2k\le n}(\beta^{-1}-\alpha^{-1})^{2k}(\beta^{-1}+\alpha^{-1})^{n-2k}
\sum_{1\le j_1<\dots<j_{2k}\le n}
\cos\kakko{2l\sum_{r=1}^{2k} (-1)^rt_{j_r}}}\\
&=2\kakko{\frac{\alpha+\beta}{2\alpha\beta}}^n
\kakko{1+%
\sum_{0<2k\le n}\kakko{\frac{\alpha-\beta}{\alpha+\beta}}^{2k}
\sum_{1\le j_1<\dots<j_{2k}\le n}
\cos\kakko{\sum_{r=1}^{2k} (-1)^rqx^2_{j_r}}}.
\end{align*}
This is the desired conclusion.
\end{proof}

If we put
\begin{align*}
\Delta_n(\vu)&\deq\begin{pmatrix}
\hA n1 & \hB1 & 0 & 0 & \dots & \hB n \\
\hB1 & \hA12 & \hB2 & 0 & \dots & 0 \\
0 & \hB2 & \hA23 & \hB3 & \dots & 0 \\
0 & 0 & \hB3 & \ddots & \ddots & \vdots\\
\vdots & \vdots & \vdots & \ddots & \ddots & \hB{n-1} \\
\hB n & 0 & 0 & \dots & \hB{n-1} & \hA{n-1}n
\end{pmatrix}\\
&=\sum_{i=1}^n\ckakko{\kakko{E^{(n)}_{ii}+E^{(n)}_{i+1,i+1}}\kakko{\frac1{1-u_i^4}-\frac12}%
+\kakko{E^{(n)}_{i,i+1}+E^{(n)}_{i+1,i}}\frac{-u_i^2}{1-u_i^4}},
\end{align*}
then
\begin{align}
\E(u_1,\dots,u_n;x_1,\dots,x_n)=\exp\kakko{-\vx\Delta_n(\vu)\vx'}
\end{align}
and
\begin{align}\label{eq:Vn}
\det\Delta_n(\vu)=\frac{(1-u_1^2\dots u_n^2)^2}{(1-u_1^4)\dots(1-u_n^4)}
\end{align}
(see \cite[Theorem A.2]{IW2005b}).
Here we denote by $E^{(n)}_{ij}$ is the matrix unit of size $n$,
and we assume that the indices of $E^{(n)}_{ij}$ are understood modulo $n$,
i.e. $E^{(n)}_{0,j}=E^{(n)}_{n,j}$, $E^{(n)}_{n+1,j}=E^{(n)}_{1,j}$, etc.
The prime $'$ indicates the matrix transpose.
Notice that $\Delta_n(\vu)$ is real symmetric and positive definite for any $\vu\in(0,1)^n$.
For $\{j_1,j_2,\dots,j_{2k}\}\subset[n]=\{1,2,\dots,n\}$, we also put
\begin{align*}
\Xi_n(j_1,\dots,j_{2k})\deq\I\sum_{r=1}^{2k}(-1)^rE^{(n)}_{j_r,j_r}.
\end{align*}
Since
\begin{align*}
\sum_{r=1}^{2k}(-1)^rx_{j_r}^2=\vx\Xi_n(j_1,\dots,j_{2k})\vx'
\end{align*}
and
\begin{align*}
\cos\kakko{q\sum_{r=1}^{2k}(-1)^rx_{j_r}^2}
=\frac12\ckakko{\exp\kakko{\I q\sum_{r=1}^{2k}(-1)^rx_{j_r}^2}+\exp\kakko{-\I q\sum_{r=1}^{2k}(-1)^rx_{j_r}^2}},
\end{align*}
we have
\begin{equation}\label{eq:Ecos}
\begin{split}
&\eqsp\E(u_1,\dots,u_n;x_1,\dots,x_n)\cos\kakko{q\sum_{r=1}^{2k}(-1)^rx_{j_r}^2}\\
&=\frac12\exp\kakko{-\vx\kakko{\Delta_n(\vu)+q\Xi_n(j_1,\dots,j_{2k})}\vx'}%
+\frac12\exp\kakko{-\vx\kakko{\Delta_n(\vu)-q\Xi_n(j_1,\dots,j_{2k})}\vx'}.
\end{split}
\end{equation}
As in \cite[Lemma A.1]{IW2005b}, one prove the
\begin{lem}\label{lem:reality}
The determinant
\begin{align}\label{eq:evenness}
\det\kakko{\Delta_n(\vu)+q\Xi_n(j_1,\dots,j_{2k})}
\end{align}
is even in $q$.
In particular, it follows that this determinant is real-valued for each $\vu\in(0,1)^n$ and $q>0$.
\qed
\end{lem}

Denote by $\cyclic m$ the cyclic subgroup of the symmetric group $\sym m$ of degree $m$
generated by the cyclic permutation $(1,2,\dots,m)\in\sym{m}$.
By Lemma \ref{lem:reality}, it follows that
\begin{align*}
\det\kakko{\Delta_n(\vu)+q\Xi_n(j_1,\dots,j_{2k})}
=\det\kakko{\Delta_n(\vu)+q\Xi_n(j_{\sigma(1)},\dots,j_{\sigma(2k)})}
\end{align*}
for any $\sigma\in \cyclic{2k}$
since $\Xi_n(j_{\sigma(1)},\dots,j_{\sigma(2k)})=\sgn(\sigma)\Xi_n(j_1,\dots,j_{2k})$.

%-------------------------------------------------------------------------
\subsection{LDU decomposition and certain positivity}
%-------------------------------------------------------------------------

Denote by $\symx n$ the set of $n$ by $n$ \emph{complex} symmetric matrices such that all principal minors are invertible,
and by $\Sym_n^+(\R)$ the set of $n$ by $n$ \emph{positive real} symmetric matrices.
Notice that $\Delta_n(\vu)\in\Sym_n^+(\R)$ for any $\vu\in(0,1)^n$.

\begin{lem}[LDU decomposition]\label{lem:LDU}
Let $n$ be a positive integer.
For any $A\in\symx n$,
there exists a lower unitriangular matrix $L$ and a diagonal matrix $D$ such that $A=LDL'$.
Moreover, $D$ is given by
\begin{align*}
D=\diag(d_1,d_2/d_1,d_3/d_2,\dots,d_n/d_{n-1}),
\end{align*}
where $d_k$ denotes the $k$-th principal minor determinant of $A$.
\end{lem}

\begin{proof}
Let us prove by induction on $n$.
The assertion is clear if $n=1$.
Suppose that the assertion is true for $n-1$.
Take $A\in\symx n$ and write
\begin{align*}
A=\begin{pmatrix} A_0 & \va \\ \va' & \alpha \end{pmatrix}
\end{align*}
with $A_0\in\symx{n-1}$, $\va\in\C^{n-1}$ and $\alpha\in\C$.
By the induction hypothesis,
there exist lower unitriangular matrix $L_0$ and diagonal matrix $D_0$ of size $n-1$ such that $A_0=L_0D_0L_0'$.
Put
\begin{align*}
L=\begin{pmatrix} L_0 & \vo \\ \vv' & 1 \end{pmatrix},\quad
D=\begin{pmatrix} D_0 & \vo \\ \vo' & d \end{pmatrix},
\end{align*}
where $\vv=(L_0D_0)^{-1}\va$ and $d=\alpha-\va'A_0^{-1}\va$
(notice that $(L_0D_0)^{-1}$ and $A_0^{-1}$ exist by the induction hypothesis)
and $\vo\in\C^{n-1}$ represents the zero vector.
Then it is straightforward to check that $A=LDL'$.
This prove the first assertion of the lemma.
The second assertion is obvious by the construction of $D$ above.
\end{proof}

\begin{lem}\label{lem:positivity_lemma}
Let $T\in\Sym_n^+(\R)$ and $D$ be a real diagonal matrix of size $n$.
Denote by $d_m$ the principal $k$-minor determinant of $T+\I D$.
Then it follows that $\Re\kakko{d_{m+1}\conj{d_m}}>0$ for $m=1,2,\dots,n-1$.
\end{lem}

\begin{proof}
Clearly, it is enough to prove the positivity of $\Re\kakko{d_{m+1}\conj{d_m}}$ with $n=m+1$.
Write $T$ and $D$ as
\begin{align*}
T=\begin{pmatrix} A & \va \\ \va' & \alpha \end{pmatrix},\quad
D=\begin{pmatrix} U & \vo \\ \vo' & u \end{pmatrix}
\end{align*}
with $A\in\Sym_m^+(\R)$, $\va\in\R^m$, $\alpha\in\R$, $u\in\R$ and a real diagonal matrix $U$ of size $m$.
Here $\vo\in\R^m$ is the zero vector.
Since $T$ is positive, we must have $0<\va'A^{-1}\va<\alpha$.
Put $B=\sqrt A\in\Sym_m^+(\R)$, $X=B^{-1}UB^{-1}\in\Sym_m(\R)$ and $\vb=B^{-1}\va$.
Then we have
\begin{align*}
d_{m+1}\conj{d_m}
&=\det\begin{pmatrix} A+\I U & \va \\ \va' & \alpha+\I u \end{pmatrix}
\begin{pmatrix} A-\I U & \vo \\ \vo & 1 \end{pmatrix}\\
&=\det\begin{pmatrix} (A+\I U)(A-\I U) & \va \\ \va'(A-\I U) & \alpha+\I u \end{pmatrix}\\
&=\abs{\det(A+\I U)}^2\kakko{\alpha+\I u-\va'(A+\I U)^{-1}\va}\\
&=\det B^4\det(I+X^2)\kakko{\alpha+\I u-\vb'\kakko{I+\I X}^{-1}\vb}.
\end{align*}
Since
\begin{align*}
\kakko{I+\I X}^{-1}=(I+X^2)^{-1}-\I X(I+X^2)^{-1},
\end{align*}
it follows that
\begin{align*}
\Re\kakko{\vb'\kakko{I+\I X}^{-1}\vb}=\vb'(I+X^2)^{-1}\vb\le\vb'\vb=\va'A^{-1}\va<\alpha
\end{align*}
or
\begin{align*}
\Re\kakko{\alpha+\I u-\vb'\kakko{I+\I X}^{-1}\vb}>0.
\end{align*}
Thus we have $\Re(d_{m+1}\conj{d_m})>0$ as desired.
\end{proof}

%=========================================================================
\section{Special values}
%=========================================================================

We recall the well-known fact.
\begin{lem}[Gaussian integral]\label{lem:gaussian}
For any $a,b\in\C$ with $\Re a>0$, it follows that
\begin{align*}
\int_{\R}\exp(-a(x-b)^2)dx=\sqrt{\frac{\pi}{a}}.
\end{align*}
Here $\sqrt a$ is chosen as $\Re\sqrt a>0$.
\qed
\end{lem}

By Lemma \ref{lem:LDU}, $A\in\symx n$ is decomposed as $A=LDL'$
with a certain lower unitriangular matrix $L$ and a diagonal matrix $D=\diag(d_1,d_2/d_1,\dots,d_n/d_{n-1})$,
where $d_k$ is the $k$-th principal minor determinant of $A$.
If all entries of $D$ have positive real parts, then it follows from Lemma \ref{lem:gaussian} that
\begin{align}\label{eq:gaussian_multidim}
\int_{\R^n}\exp(-\vx A\vx')d\vx=\frac{\pi^{n/2}}{\sqrt{\det A}}.
\end{align}
Here $\vx=(x_1,\dots,x_n)$ and $d\vx=dx_1\dotsb dx_n$.

Now the matrix $\Delta_n(\vu)+q\Xi_n(\vj)$ belongs to $\symx n$ for any $\vu\in(0,1)^n$.
Denote by $d_k=d_k(n,\vu,q,\vj)$ the $k$-th principal minor determinant of $\Delta_n(\vu)+q\Xi_n(\vj)$, and put $d_0=1$.
It then follows from Lemma \ref{lem:positivity_lemma} that $\Re(d_k/d_{k-1})>0$ for $k=1,2,\dots,n$.
Consequently, in view of \eqref{eq:Ecos}, \eqref{eq:evenness}, \eqref{eq:gaussian_multidim} and Lemma \ref{lem:expansion_of_B},
we can calculate $\F(u_1,\dots,u_n)$ as
\begin{align*}
&\F(u_1,\dots,u_n)=2\sqrt{\pi^n}\kakko{\frac{\alpha+\beta}{2\alpha\beta}}^{\!\!n}\\
&\times\ckakko{\frac1{\sqrt{\det\Delta_n(\vu)}}%
+\sum_{0<2k\le n}\kakko{\frac{\alpha-\beta}{\alpha+\beta}}^{\!\!2k}%
\sum_{1\le j_1<j_2<\dots<j_{2k}\le n}
\frac1{\sqrt{\det\kakko{\Delta_n(\vu)+q\Xi_n(j_1,\dots,j_{2k})}}}}.
\end{align*}
We also notice that
\begin{align*}
\int_{[0,1]^n}\frac{2^nd\vu}{1-u_1^2u_2^2\dotsb u_n^2}=\zeta(n,1/2)
\end{align*}
for $n\ge2$.
From these equations together with \eqref{eq:special_value_primitive} and \eqref{eq:Vn}, we now obtain the
\begin{thm}\label{thm:specialvalues}
For each positive integer $n\ge2$, it follows that
\begin{equation}\label{eq:specialvalues}
\begin{split}
\zeta_Q(n)&=2\kakko{\frac{\alpha+\beta}{2\sqrt{\alpha\beta(\alpha\beta-1)}}}^{\!\!n}
\Biggl(\zeta(n,1/2)+\sum_{0<2k\le n}\kakko{\frac{\alpha-\beta}{\alpha+\beta}}^{\!\!2k}R_{n,k}(q)\Biggr).
\end{split}
\end{equation}
Here $R_{n,k}(q)$ is given by
\begin{align*}
R_{n,k}(q)&=
\sum_{1\le j_1<j_2<\dots<j_{2k}\le n}\int_{[0,1]^n}
\frac{2^ndu_1\dots du_n}{\sqrt{\den n(\vu;q;j_1,\dots,j_{2k})}},\\
\den n(\vu;q;j_1,\dots,j_{2k})&=\det\kakko{\Delta_n(\vu)+q\Xi_n(j_1,\dots,j_{2k})}\prod_{j=1}^n(1-u_j^4).
\end{align*}
\qed
\end{thm}

\begin{rem}\label{rem:degenerate_case}
If $\alpha=\beta$, then we have $\zeta_Q(n)=2(\alpha^2-1)^{-n/2}\zeta(n,1/2)$,
which is a special case of the fact that
$\zeta_Q(s)=2(\alpha^2-1)^{-s/2}\zeta(s,1/2)$ for $\alpha=\beta$.
In fact, when $\alpha$ and $\beta$ are equal,
we can show that $Q\cong\sqrt{\alpha^2-1}\kakko{-\frac12\frac{d^2}{dx^2}+\frac12x^2}I$ (see \cite{PW2002}).
\end{rem}

\begin{rem}
Put $\omega=\sqrt{\alpha/\beta}$.
Consider the limit $\e\to+0$ under the condition that $\omega$ is kept constant.
Then we see from \eqref{eq:specialvalues} that
\begin{align}\label{eq:limit_with_constant_omega}
\frac{\zeta_Q(n)}{\e^n}\to(\omega^n+\omega^{-n})\zeta(n,1/2)\qquad(\e\to+0)
\end{align}
if $n$ is a positive integer greater than one.
Further, for each fixed $s>1$, one can prove that
$\zeta_Q(s)/\e^s$ is analytic with respect to $\e$ near the origin 
and its constant term is $(\omega^s+\omega^{-s})\zeta(s,1/2)$.
The detailed discussion on this result will be given elsewhere \cite{KP2009}.
\end{rem}

We give an expansion of the determinants $\den n(\vu;q;j_1,\dots,j_{2k})$ appearing in \eqref{eq:specialvalues}.
For $\vj=\{j_1,j_2,\dots,j_r\}\subset[n]$ with $r>0$ and $j_1<j_2<\dots<j_r$, define
\begin{align*}
C_{n}(\vu;\vj)=\prod_{i=1}^r(1-u_{j_i}^4u_{j_i+1}^4\dots u_{j_{i+1}-1}^4).
\end{align*}
We also define $C_{n}(\vu;\emptyset)=(1-u_1^2u_2^2\dots u_n^2)^2$.
Here we regard that $j_{r+1}=n+j_1$ and $u_{i+n}=u_i$.
For instance, if $n=9$ and $\vj=\{3,6,8\}$, then
\begin{align*}
C_{9}(\vu;\vj)=(1-u_3^4u_4^4u_5^4)(1-u_6^4u_7^4)(1-u_8^4u_9^4u_1^4u_2^4).
\end{align*}

\begin{lem}
For a given subset $\vj=\{j_1,j_2,\dots,j_r\}\subset[n]$ with $j_1<j_2<\dots<j_r$, it follows that
\begin{align}
\den n(\vu;q;\vj)=\sum_{d\ge0}(-q^2)^{d}\,\den{n,d}(\vu;\vj)
\end{align}
with
\begin{align}
\den{n,d}(\vu;\vj)\deq\sum_{\substack{S\subset[2k]\\\card{S}=2d}}(-1)^{\esum S}C_{n}(\vu;\vj(S)).
\end{align}
Here $\esum S\deq\sum_{s\in S}s$ is the sum of the elements in $S$ and
$\vj(S)\deq\{j_{s_1},\dots,j_{s_l}\}$ if $S=\{s_1,\dots,s_l\}$ with $s_1<\dots<s_l$.
\end{lem}

\begin{proof}
Denote by $\vd_i$ the $i$-th column vector of $\Delta_n(\vu)$.
We also denote by $\{\ve_i\}_{i=1}^n$ the standard basis of $\C^n$.
By the multilinearity of a determinant, we readily get
\begin{align*}
&\eqsp\det(\Delta_n(\vu)+q\Xi_n(\vj))\\
&=\det\Delta_n(\vu)+\sum_{r=1}^{2k}(\I q)^r
\sum_{1\le s_1<\dots<s_r\le 2k}(-1)^{s_1+\dots+s_r}\det(\vd_1,\dots,\ve_{j_{s_1}},\dots,\ve_{j_{s_r}},\dots,\vd_n).
\end{align*}
The determinant $\det(\vd_1,\dots,\ve_{j_{s_1}},\dots,\ve_{j_{s_r}},\dots,\vd_n)$ is a product of
$r$ tridiagonal determinants
\begin{align*}
D_i=
\begin{vmatrix}
d_{a(i)+1,a(i)+1} & d_{a(i)+1,a(i)+2} & \dots & d_{a(i)+1,a(i+1)-1} \\
d_{a(i)+2,a(i)+1} & d_{a(i)+2,a(i)+2} & \dots & d_{a(i)+2,a(i+1)-1} \\
\vdots & \vdots & \ddots & \vdots \\
d_{a(i+1)-1,a(i)+1} & d_{a(i+1)-1,a(i)+2} & \dots & d_{a(i+1)-1,a(i+1)-1}
\end{vmatrix},
\end{align*}
where $a(i)=j_{s_i}$, $d_{ij}$ is the $(i,j)$-entry of $\Delta_n(\vu)$,
and the indices are understood modulo $n$.
If $a(i+1)=a(i)+1$, then we understand that $D_i=1$.
It is easy to see that
\begin{align*}
D_i=\frac{1-u_{a(i)}^4u_{a(i)+1}^4\dotsb u_{a(i+1)-1}^4}{(1-u_{a(i)}^4)(1-u_{a(i)+1}^4)\dotsb(1-u_{a(i+1)-1}^4)}.
\end{align*}
Hence we have
\begin{align*}
\den n(\vu;q;\vj)=\sum_{S\subset[2k]}(-1)^{\esum S}(\I q)^{\card{S}}C_{n}(\vu;\vj(S)).
\end{align*}
Since $\den n(\vu;q;\vj)$ is real-valued by Lemma \ref{lem:reality},
we have the conclusion by taking the real parts.
\end{proof}

%=========================================================================
\section{Examples}
%=========================================================================

%-------------------------------------------------------------------------
\subsection{$\den{n,d}(\vu;\vj)$ and $R_{n,k}(q)$}
%-------------------------------------------------------------------------

We give several examples of $\den{n,d}(\vu;\vj)$.
For convenience, we prepare some notation for abbreviation.
Let us put
\begin{align*}
\V n(\vu)\deq(1-u_1^2\dots u_n^2)^2,\qquad
\U\vt(\vu)\deq\prod_{i=1}^m\biggl(1-\prod_{j=1}^{t_i}u^4_{j+\sum_{k<i}t_k}\biggr)
\end{align*}
for a positive integer $n$ and a sequence $\vt=(t_1,\dots,t_m)\in\T m(n)$, where
\begin{align*}
\T m(n)\deq\ckakko{(t_1,\dots,t_m)\in[n]^4\,;\,t_1+\dotsb+t_m=n}.
\end{align*}
For instance, if $\vt=(2,3,2,1)\in\T4(8)$, then
\begin{align*}
\U\vt(u_1,\dots,u_8)=(1-u_1^4u_2^4)(1-u_3^4u_4^4u_5^4)(1-u_6^4u_7^4)(1-u_8^4).
\end{align*}
Notice that $\den{n,0}(\vu;\vj)=\V n(\vu)$ for any $\vj$.

\begin{ex}\label{ex:Rn1}
For $\vj=\{i,j\}\subset[n]$ with $i<j$, we have
\begin{align*}
\den{n,1}(\vu;\vj)=(-1)^{1+2}C_{n}(\vu;\vj)=-\U{(r,n-r)}(u_i,u_{i+1},\dots,u_{i-1}),
\end{align*}
where $r=j-i$.
This fact immediately implies that $R_{n,1}(q)$ in \eqref{eq:specialvalues} is given by
\begin{align*}
R_{n,1}(q)
&=\frac n2\sum_{r=1}^{n-1}\int_{[0,1]^n}\frac{2^nd\vu}{\sqrt{\V n(\vu)+q^2\U{(r,n-r)}(\vu)}}\\
&=\sum_{0<2r\le n}\frac n{1+\delta_{2r,n}}\int_{[0,1]^n}\frac{2^nd\vu}{\sqrt{\V n(\vu)+q^2\U{(n-r,r)}(\vu)}}.
\end{align*}
\end{ex}

\begin{ex}\label{ex:top_term}
For $\vj\subset[n]$ with $\card\vj=2k$, it follows in general that
\begin{align*}
\den{n,k}(\vu;\vj)=(-1)^{k}C_{n}(\vu;\vj)
\end{align*}
since $\esum{[2k]}=k(2k-1)\equiv k\pmod2$.
\end{ex}

\begin{ex}\label{ex:Rn2}
For $\vj=\{j_1,j_2,j_3,j_4\}\subset[n]$ with $j_1<j_2<j_3<j_4$, we have
\begin{align*}
\den{n,1}(\vu;\vj)
&=(-1)^{1+2}C_{n}(\vu;j_1,j_2)+(-1)^{1+3}C_{n}(\vu;j_1,j_3)+(-1)^{1+4}C_{n}(\vu;j_1,j_4)\\
&\eqsp+(-1)^{2+3}C_{n}(\vu;j_2,j_3)+(-1)^{2+4}C_{n}(\vu;j_2,j_4)+(-1)^{3+4}C_{n}(\vu;j_3,j_4)\\
&=-(1-u_{j_1}^4\dots u_{j_2-1}^4)(1-u_{j_3}^4\dots u_{j_4-1}^4)(1-u_{j_2}^4\dots u_{j_3-1}^4)(1-u_{j_4}^4\dots u_{j_1-1}^4)\\
&\eqsp-(1-u_{j_1}^4\dots u_{j_2-1}^4u_{j_3}^4\dots u_{j_4-1}^4)(1-u_{j_2}^4\dots u_{j_3-1}^4u_{j_4}^4\dots u_{j_1-1}^4).
\end{align*}
By Example \ref{ex:top_term}, we also see that
\begin{align*}
\den{n,2}(\vu;\vj)&=C_n(\vu;\vj)\\
&=(1-u_{j_1}^4\dots u_{j_2-1}^4)(1-u_{j_2}^4\dots u_{j_3-1}^4)(1-u_{j_3}^4\dots u_{j_4-1}^4)(1-u_{j_4}^4\dots u_{j_1-1}^4).
\end{align*}
Thus we have
\begin{align*}
&\eqsp\det(\Delta_n(\vu)+q\Xi_n(j_1,j_2,j_3,j_4))\prod_{i=1}^n(1-u_i^4)\\
&=\V n(\vu)+(q^2+q^4)(1-u_{j_1}^4\dots u_{j_2-1}^4)(1-u_{j_3}^4\dots u_{j_4-1}^4)(1-u_{j_2}^4\dots u_{j_3-1}^4)(1-u_{j_4}^4\dots u_{j_1-1}^4)\\
&\eqsp+q^2(1-u_{j_1}^4\dots u_{j_2-1}^4u_{j_3}^4\dots u_{j_4-1}^4)(1-u_{j_2}^4\dots u_{j_3-1}^4u_{j_4}^4\dots u_{j_1-1}^4).
\end{align*}
If we take $(t_1,t_2,t_3,t_4)\in\T4(n)$ such that $j_{i+1}\equiv j_i+t_i\pmod{n}$ ($i=1,2,3,4$; $j_5=j_1$),
then it follows that
\begin{multline*}
\int_{[0,1]^n}\frac{d\vu}{\sqrt{\det(\Delta_n(\vu)+q\Xi_n(j_1,j_2,j_3,j_4))\prod_{i=1}^n(1-u_i^4)}}\\
=\int_{[0,1]^n}\frac{d\vu}{\sqrt{\V n(\vu)+q^2\U{(t_1+t_3,t_2+t_4)}(\vu)+(q^2+q^4)\U{(t_1,t_2,t_3,t_4)}(\vu)}}.
\end{multline*}
The cyclic group $\cyclic 4$ of order $4$ naturally acts on $\T4(n)$ by
\begin{align*}
\sigma.(t_1,t_2,t_3,t_4)\deq(t_{\sigma(1)},t_{\sigma(2)},t_{\sigma(3)},t_{\sigma(4)}) \quad (\sigma\in\cyclic4).
\end{align*}
Notice that the integral above is $\cyclic 4$-invariant.
For a given $\vt=(t_1,t_2,t_3,t_4)\in\T4(n)$,
the number of subsets $\vj=\{j_1,j_2,j_3,j_4\}$ in $[n]$
satisfying the condition $j_{i+1}\equiv j_i+t_i\pmod{n}$ is equal to $n/\card{C_4(\vt)}$,
where $\cyclic4(\vt)$ denotes the stabilizer of $\vt$ in $\cyclic4$.
Consequently,
\begin{align*}
R_{n,2}(q)
&=\sum_{\vt\in\T4(n)/\cyclic4}\frac n{\card{\cyclic4(\vt)}}
\int_{[0,1]^n}\frac{2^nd\vu}{\sqrt{\V n(\vu)+q^2\U{(t_1+t_3,t_2+t_4)}(\vu)+(q^2+q^4)\U{(t_1,t_2,t_3,t_4)}(\vu)}}\\
&=\frac n4\sum_{\vt\in\T4(n)}
\int_{[0,1]^n}\frac{2^nd\vu}{\sqrt{\V n(\vu)+q^2\U{(t_1+t_3,t_2+t_4)}(\vu)+(q^2+q^4)\U{(t_1,t_2,t_3,t_4)}(\vu)}},
\end{align*}
where $\vt=(t_1,t_2,t_3,t_4)$.
Similarly, the result in Example \ref{ex:Rn1} can be also rewritten as
\begin{align*}
R_{n,1}(q)
&=\frac n2\sum_{\vt\in\T2(n)}\int_{[0,1]^n}\frac{2^nd\vu}{\sqrt{\V n(\vu)+q^2\U{(t_1,t_2)}(\vu)}}\\
&=\sum_{\vt\in\T2(n)/\cyclic2}\frac n{\card{\cyclic2(\vt)}}\int_{[0,1]^n}\frac{2^nd\vu}{\sqrt{\V n(\vu)+q^2\U{(t_1,t_2)}(\vu)}}.
\end{align*}
\end{ex}

%-------------------------------------------------------------------------
\subsection{Several special values}
%-------------------------------------------------------------------------

Using Theorem \ref{thm:specialvalues} and the formulas for $R_{n,1}(q)$ and $R_{n,2}(q)$ given in the previous example,
we show several examples of the special values of $\zeta_Q(s)$.

\begin{ex}[Ichinose-Wakayama's result]
The values $\zeta_Q(2)$ and $\zeta_Q(3)$ are given by
\begin{align*}
\zeta_Q(2)&=2\kakko{\frac{\alpha+\beta}{2\sqrt{\alpha\beta(\alpha\beta-1)}}}^{\!\!2}
\Biggl(\zeta(2,1/2)+\kakko{\frac{\alpha-\beta}{\alpha+\beta}}^{\!\!2}R_{2,1}(q)\Biggr),\\
\zeta_Q(3)&=2\kakko{\frac{\alpha+\beta}{2\sqrt{\alpha\beta(\alpha\beta-1)}}}^{\!\!3}
\Biggl(\zeta(3,1/2)+\kakko{\frac{\alpha-\beta}{\alpha+\beta}}^{\!\!2}R_{3,1}(q)\Biggr)
\end{align*}
with
\begin{align*}
R_{2,1}(q)&=\int_{[0,1]^2}\frac{4du_1du_2}{\sqrt{\V2(\vu)+q^2\U{(1,1)}(\vu)}}
=\int_{[0,1]^2}\frac{4du_1du_2}{\sqrt{(1-u_1^2u_2^2)^2+q^2(1-u_1^4)(1-u_2^4)}},\\
R_{3,1}(q)&=3\int_{[0,1]^3}\frac{8du_1du_2du_3}{\sqrt{\V3(\vu)+q^2\U{(2,1)}(\vu)}}
=3\int_{[0,1]^3}\frac{8du_1du_2du_3}{\sqrt{(1-u_1^2u_2^2u_3^2)^2+q^2(1-u_1^4)(1-u_2^4u_3^4)}}.
\end{align*}
This recovers the result obtained by Ichinose and Wakayama \cite{IW2005b}.
\end{ex}

\begin{ex}
The values $\zeta_Q(4)$ and $\zeta_Q(5)$ are given by
\begin{align*}
\zeta_Q(4)&=2\kakko{\frac{\alpha+\beta}{2\sqrt{\alpha\beta(\alpha\beta-1)}}}^{\!\!4}
\Biggl(\zeta(4,1/2)+\kakko{\frac{\alpha-\beta}{\alpha+\beta}}^{\!\!2}R_{4,1}(q)
+\kakko{\frac{\alpha-\beta}{\alpha+\beta}}^{\!\!4}R_{4,2}(q)\Biggr),\\
\zeta_Q(5)&=2\kakko{\frac{\alpha+\beta}{2\sqrt{\alpha\beta(\alpha\beta-1)}}}^{\!\!5}
\Biggl(\zeta(5,1/2)+\kakko{\frac{\alpha-\beta}{\alpha+\beta}}^{\!\!2}R_{5,1}(q)
+\kakko{\frac{\alpha-\beta}{\alpha+\beta}}^{\!\!4}R_{5,2}(q)\Biggr)
\end{align*}
with
\begin{align*}
R_{4,1}(q)&=4\int_{[0,1]^4}\frac{16d\vu}{\sqrt{\V4(\vu)+q^2\U{(3,1)}(\vu)}}
+2\int_{[0,1]^4}\frac{16d\vu}{\sqrt{\V4(\vu)+q^2\U{(2,2)}(\vu)}},\\
R_{5,1}(q)&=5\int_{[0,1]^5}\frac{32d\vu}{\sqrt{\V5(\vu)+q^2\U{(4,1)}(\vu)}}
+5\int_{[0,1]^5}\frac{32d\vu}{\sqrt{\V5(\vu)+q^2\U{(3,2)}(\vu)}}
\end{align*}
and
\begin{align*}
R_{4,2}(q)&=\int_{[0,1]^4}\frac{16d\vu}{\sqrt{\V4(\vu)+q^2\U{(2,2)}(\vu)+(q^2+q^4)\U{(1,1,1,1)}(\vu)}},\\
R_{5,2}(q)&=5\int_{[0,1]^5}\frac{32d\vu}{\sqrt{\V5(\vu)+q^2\U{(3,2)}(\vu)+(q^2+q^4)\U{(2,1,1,1)}(\vu)}}.
\end{align*}
\end{ex}

\begin{rem}[Ap\'ery-like numbers and the Heun differential equation]\label{rem:Apery-Heun}
If we define the numbers $J_2(m)$ ($m\ge0$) by the expansion
\begin{align*}
R_{2,1}(q)=\sum_{m=0}^\infty \binom{-1/2}m J_2(m)q^{2m},
\end{align*}
then they satisfy the three-term recurrence relation
\begin{align*}
4m^2J_2(m)-(8m^2-8m+3)J_2(m-1)+4(m-1)^2J_2(m-2)=0\quad(m\ge1).
\end{align*}
This implies that the generating function $w_2(t)=\sum_{m=0}^\infty J_2(m)t^m$ satisfies
\begin{align}\label{eq:Heun}
\ckakko{t(1-t)^2\frac{d^2}{dt^2}+(1-3t)(1-t)\frac{d}{dt}+t-\frac34}w_2(t)=0,
\end{align}
which is a singly confluent Heun differential equation \cite{IW2005b}.
Fortunately, this equation can be reduced to the Gaussian hypergeometric defferential equation
by a suitable change of variable and solved as follows \cite{O}:
\begin{align*}
w_2(t)=\frac{3\zeta(2)}{1-t}{}_2F_1\kakko{\frac12,\frac12;1;\frac t{t-1}},
\end{align*}
from which we obtain
\begin{align*}
R_{2,1}(q)=3\zeta(2){}_2F_1\kakko{\frac14,\frac34;1;-q^2}^{\!2}.
\end{align*}
Thus we have the following formulas for $\zeta_Q(2)$ \cite{IW2005b,O}:
\begin{align*}
\zeta_Q(2)
&=\kakko{\frac{\pi(\alpha+\beta)}{2\sqrt{\alpha\beta(\alpha\beta-1)}}}^{\!\!2}
\Biggl(1+\frac1{2\pi\I}\kakko{\frac{\alpha-\beta}{\alpha+\beta}}^{\!\!2}\int_{\abs z=r}\frac{u(z)}{z(1+q^2z)^{1/2}}dz\Biggr)\\
&=\kakko{\frac{\pi(\alpha+\beta)}{2\sqrt{\alpha\beta(\alpha\beta-1)}}}^{\!\!2}
\Biggl(1+\kakko{\frac{\alpha-\beta}{\alpha+\beta}}^{\!\!2}{}_2F_1\kakko{\frac14,\frac34;1;-q^2}^{\!2}\Biggr),
\end{align*}
where $u(z)=w_2(z)/3\zeta(2)$ is a normalized (unique) holomorphic solution of \eqref{eq:Heun} in $\abs z<1$ and $q^2<r<1$.
We also have similar formulas for $\zeta_Q(3)$ \cite{IW2005b,KW2006}.
Furthermore, if we define $J_n(m)$ by the expansion of the integral
\begin{align*}
\int_{[0,1]^n}\frac{2^nd\vu}{\sqrt{\V n(\vu)+q^2\U{(n-1,1)}(\vu)}}=\sum_{m=0}^\infty \binom{-1/2}m J_n(m)q^{2m}
\end{align*}
appearing in $R_{n,1}(q)$, we also have the recurrence relation
\begin{align*}
4m^2J_n(m)-(8m^2-8m+3)J_n(m-1)+4(m-1)^2J_n(m-2)=J_{n-2}(m)\quad(m\ge1,\ n\ge4)
\end{align*}
and the differential equation
\begin{align*}
\ckakko{t(1-t)^2\frac{d^2}{dt^2}+(1-3t)(1-t)\frac{d}{dt}+t-\frac34}w_n(t)=\frac{w_{n-2}(t)-w_{n-2}(0)}{4t}\quad(n\ge4)
\end{align*}
for $w_n(t)=\sum_{m=0}^\infty J_n(m)t^m$ \cite{K2009}.
However, of course, these results are still insufficient to obtain simpler and/or more explicit formula for $\zeta_Q(n)$ for $n\ge4$.
\end{rem}

\begin{ackn}
The author would like to thank Max-Planck-Institut f\"ur Mathematik for the support and hospitality.
%during his stay from December 2008 to March 2009.
%Indeed, the first draft of this paper was completed in that stay.
The author is also partially supported by Grand-in-Aid for Young Scientists (B) No. 20740021.
\end{ackn}

%=========================================================================
%	References
%=========================================================================

\bigskip
\noindent
Department of Mathematical Sciences, University of the Ryukyus\\
1 Senbaru, Nishihara-cho, Okinawa 903-0213 Japan

\smallskip
\noindent
\texttt{kimoto@math.u-ryukyu.ac.jp}

\smallskip
\noindent
Max-Planck-Institut f\"ur Mathematik\\
Vivatsgasse 7, 53111 Bonn, Germany

\smallskip
\noindent
\texttt{kimoto@mpim-bonn.mpg.de}

\end{document}